\documentclass[12pt]{article}
\usepackage{amssymb}
\def\hang{\hangindent\parindent}
\def\tex#1{\indent\llap{[#1]\enspace}\ignorespaces}
\def\re{\par\hang\tex}

\def\a{\alpha}
\def\b{\beta}
\def\d{\delta}

\def\WW{{\cal W}}
\def\g{\gamma}
\def\G{\Gamma}
\def\l{\lambda}

\def\sc{\scriptstyle}
\def\ssc{\scriptscriptstyle}
\def\dis{\displaystyle}
\def\cl{\centerline}

\def\ol{\overline}

\def\sptl{\check\ptl}

\def\bs{\backslash}

\def\vs{\vspace*}

\def\AA{{\cal A}}
\def\DD{{\cal D}}

\def\ni{\noindent}
\def\ptl{\partial}

\def\Z{\mathbb{Z}{\ssc\,}}

\def\F{\mathbb{F}{\ssc\,}}
\begin{document}
\cl{{\Large \bf
 2-Cocycles of the Lie superalgebras of Weyl type}\footnote{Supported by NSF grant 10171064 of China}}
\vs{6pt}

\cl{(appeared in {\it Comm. Algebra}, {\bf33} (2005), 2991--3007)}
\vskip6pt
 \cl{ Guang'ai
Song and Yucai Su} \cl{\small Department of Mathematics, Shanghai
Jiao Tong University,
 Shanghai, 200240, \vs{-4pt}China}
\cl{\small E-mail: gasong@sjtu.edu.cn, ycsu@sjtu.edu.cn}
\vs{6pt}

\noindent{\small{\bf Abstract.} In a paper by Su
and Zhao, the Lie algebra $\AA[\DD]=\AA\bigotimes \F[\DD]$ of Weyl type was
defined and studied, where $\AA$ is a commutative associative
algebra with an identity element over a field $\F$ of arbitrary
characteristic, and $\F[\DD]$ is the polynomial algebra of a commutative
derivation subalgebra $\DD$ of $\AA$. The 2-cocycles of a class of $\AA[\DD]$ were
determined by Su. In the present paper, we determine the 2-cocycles
of a class of Lie superalgebras of Weyl type over a field $\F$ of
characteristic 0.

}
\vs{6pt}

\cl{\bf\S1. \ Introduction}\setcounter{section}{1}\setcounter{equation}{0}
Recently there appeared a number of papers on the
structure theory of infinite dimensional
Lie (super)algebras and conformal (super)algebras (for example,
[SXZ,  SZ1--SZ3, SZZ, \mbox{X1--X5}, Z] and references cited in those papers).
Among them, a class of Lie algebras of Weyl type,
which is closely related to $\WW$-infinity algebras $\WW_{1+\infty},\,\WW_\infty$ and the general conformal algebras $gc_N$ (see, e.g., [BKV, S3, S4]), was studied
in [SZ1, SZZ, Z].
In this \mbox{paper,} we study the $2$-cohomology groups of
the Lie \mbox{superalgebras} of
Weyl type, a natural generalization of Lie algebras of Weyl type,
which are
closely related to $\WW$-infinity superalgebras and general conformal \mbox{superalgebras.}
The main result of this paper is Theorem 3.5.
Since the classification of finite dimensional simple Lie superalgebras was given in [K2], the role of Lie superalgebras has become more and more important in solving problems in the quantum field theory and string theory.
A cohomology theory of Lie superalgebras and Lie color algebras was
developed in [ScZ1], while
a general theory of cohomology of Lie conformal algebras was
established in [BKV]. The
cohomology theory of Lie algebras has played important roles in
the structure and representation theories of Lie algebras.
It is well-known that central extensions, which are determined by
$2$-cohomology groups,
are widely used in the structure theory and
the representation theory of Lie algebras (e.g., [K1, KP, S3]).
Unlike the finite dimensional
 simple Lie algebras, since the complete reducibility of modules does
not hold in super case or conformal case, one may expect that, as pointed in [BKV], the cohomology
theory of Lie superalgebras and conformal superalgebras,
is very interesting and important, and further, one may expect that it is more difficult as is already seen even for the case of finite dimensional classical simple Lie superalgebras $sl(m/n)$ [ScZ1, ScZ2].
\par
The problem of determining the cohomology groups of general conformal algebras $gc_N$ remains open (see [BKV]). Using some techniques developed in [S2] which determined $2$-cocycles of Lie algebras of Weyl type, a partial answer to the problem was given
in [S4]. One of our motivation in this paper is to develop some results and techniques
in order to be used
to determine cohomology groups of the general conformal superalgebras in the future.
\par
Now we begin with some brief definitions. Let $\F$ be a field
of characteristic $0$. For any positive integer $\ell$, an additive subgroup
$\G$ of the $\ell$-dimensional vector space $\F^\ell$ is called
{\it nondegenerate} if $\G$ contains an $\F$-basis
of $\F^\ell$. Let $\ell_1, \ell_2, \cdots, \ell_5 $ be five nonnegative
integers such that $\ell=\sum_{p=1}^5 \ell_p>0$. For convenience, we
denote
\begin{eqnarray}
\mbox{$
\ell_i ^{\prime} =\sum\limits_{p=1} ^i \ell_p,
$}
\mbox{ \ and \ }
\overline{m, n}=\{m, m+1, \cdots, n\} \ {\rm if } \ m,n\in\Z,\,m\leq n.
\end{eqnarray}
An element of $\F^\ell$ will be written as $\a=(\a_1,...,\a_\ell)$.
 Take an additive subgroup $\G$ of $\F^\ell$ such that
\begin{eqnarray}
\a=(\a_1, \cdots, \a_\ell)=(0, \cdots,0, \a_{\ell_1 + 1 },
\cdots, \a_{\ell_4^{\prime}}, 0 \cdots, 0)
\in\{0\}^{\ell_1}\times\F^{\ell_2+\ell_3+\ell_4}\times\{0\}^{\ell_5}
\end{eqnarray}
for all $\a \in \G$, and such that
$\G$ is nondegenerate as a subgroup of $\F^{\ell_2+\ell_3+\ell_4}$. Set
\begin{eqnarray}
\label{1.3}
\vec{J}={\Z_+ ^{\ell_2 ^{\prime}}}\times{\Z^{\ell_3}}\times
{\{0\}^{\ell_4}}\times{\{0, 1\}^{\ell_5}}.
\end{eqnarray} For convenience, we shall always decompose a vector $\vec k=(k_1,...,k_\ell)\in \vec J$ as
\begin{eqnarray}
\label{add-k}
&&\vec{k}=\vec{i}+\vec{j},\mbox{ where,}\\
\nonumber&&\vec{i}=(i_1,...,i_\ell)=(k_1,
\cdots, k_{\ell_3^{\prime}}, 0, \cdots, 0)\in \Z^{\ell'_3}\times\{0\}^{\ell_4+\ell_5},\\
\nonumber&&
\vec{j}=(j_1,...,j_\ell)=(0,
\cdots, 0, k_{\ell_4^{\prime} +1}, \cdots, k_\ell)\in
\{0\}^{\ell_4^{\prime}}\times \{0, 1\}^{\ell_5},
\end{eqnarray}
  and we denote
\begin{eqnarray}
&&\mbox{$|\vec{k}|=\sum\limits_{p=1} ^\ell |k_p|,$}
\ \ \ \
a_{[p]}=\stackrel{p}{(0, \cdots, 0, a, 0, \cdots, 0)},
\end{eqnarray}
where $a\in \Z $ if $p\in \overline{1, \ell_4 ^{\prime}}$,  and $a\in\{0,
1\}$ if  $p\in \overline{\ell_4 ^{\prime} +1, \ell}$.

Let
\begin{equation}\label{add-A-1}
\AA_1=\F[\G]={\rm span}\{x^\a\,|\,\a\in\G\}
\end{equation}
 be the group algebra with product
$x^\a
x^\b=x^{\a+\b}$ for $\a,\b\in\G$.
Let
\begin{equation}
\label{add-A-2}
\AA_2=\F[\vec J]=\F[t_1,t_2,\cdots,t_{\ell'_3},\,s_1,s_2,\cdots,s_{\ell_5}]
\end{equation}
be the (super)polynomial algebra with the {\it ordinary} (or {\it even}) {\it variables}
$t_p,\,p\in\ol{1,\ell'_3}$ and the
{\it Grassmannian} (or {\it odd}) {\it variables}
$s_q,\,q\in\ol{1,\ell_5}$, namely, variables $t_p,\,s_q$ satisfy
\begin{equation}
\label{add-A-3}
t_p\/ a=a\/ t_p,\ \ \
s_{q}s_{q'}=-s_{q'}s_{q}
\mbox{\ (in particular, \ $s_q^2=0$)}
\end{equation}
for all $a\in\AA_2,\,p\in\ol{1,\ell_3},\,q,q'\in\ol{1,\ell_5}.$
Let $\AA=\AA_1\otimes\AA_2$, and denote
\begin{equation}
\label{denote1}
x^{\a,\vec k}=x^\a {\sc\,}t_1^{i_1}{\sc\,}\cdots {\sc\,}t_{\ell'_3}^{i_{\ell'_3}}{\sc\,}s_1^{j_{\ell'_4+1}}{\sc\,}\cdots {\sc\,}s_{\ell_5}^{j_\ell},\ \
x^{0, 0} =1,\ \
t^{\vec{i}} = x^{0, \vec{i}},\ \
s^{\vec j} = x^{0,\vec{j}},
\end{equation}
for $\a\in\G$ and $\vec k\in\vec J$ written as in (\ref{add-k}).
Then $\AA=\F[\G\times\vec J]$ is a semi-group superalgebra with basis
$\{x^{\a, \vec{k}}\,|\,(\a, \vec{k}) \in
\G\times\vec{J}\}$, and the product
\begin{eqnarray}
\label{add1.6}
x^{\a, \vec{k}}
x^{\a^{\prime}, \vec{k^{\prime}}}
=x^{\a, \vec{i}+\vec{j}} 
x^{\a^{\prime},
\vec{i^{\prime}}+\vec{j^{\prime}}} =(-1)^{\sum_{\ell'_4<p<q\le\ell} j_q j_p
^{\prime}} x^{\a +\a^{\prime}, \vec{k}+\vec{k^{\prime}}}
\end{eqnarray}
for $\a,\a'\in\G,\,\vec k,\vec k'\in\vec J,$
where we have used the following convention.

{\bf Convention 1.1.} \
If an undefined notion appears in an expression, we regard it as zero; for instance,
$x^{\a,\vec k}=0$ if $k_p\ge2$ for some $p\in\ol{\ell'_4,\ell}$ since in this case
$\vec k\notin\vec J$.

Let $\Z_2=\Z/2\Z=\{\bar0,\bar1\}$.
Then $\AA=\AA_{\bar0}+\AA_{\bar1}$ is a $\Z_2$-graded supercommutative superalgebra with the gradation spaces given by
\begin{eqnarray}
\label{Gradation-of-A}
\AA_{\bar{0} }=
{\rm span}\{x^{\a, \vec{i}+\vec{j}}\,\,{\biggl|}\,\,|\vec{j}| \mbox{ is even}\}, \ \ \
\AA_{\bar{1} }={\rm span}\{x^{\a,\vec{i}+\vec{j}}\, \,{\biggr|}\,\,
|\vec{j}| \mbox{ is odd}\}.
\end{eqnarray}
Define the linear transformations
$\{\ptl_1 ^-,
\ptl_2^-, \cdots,
\ptl^-_{\ell_3 ^{\prime}}, \ptl^+_{\ell_3 ^{\prime} +1}, \cdots,
\ptl_{\ell_4 ^{\prime}} ^+ , \ptl_{\ell_4^{\prime} +1}^-
, \cdots, \ptl_\ell ^- \}$ on $\AA$ by
\begin{equation}
\ptl_p ^- (x^{\a, \vec{k}}) =k_p x^{\a, \vec{k-1_{[p]}}},
\ \ \  \ptl_{p^{\prime}} ^+ (x^{\a,
\vec{k}})=\a_{p^{\prime}} x^{\a, \vec{k}},
\ \ \
\ptl_q ^-(x^{\a, \vec{k}}) =(-1)^{\sum_{r=\ell_4 ^{\prime}
+1} ^{q-1} k_r} k_q x^{\a, \vec{k } -1_{[q]}},
\end{equation}
for $p\in \overline{1, \ell_3 ^{\prime}}$, $ p^{\prime} \in
\overline{\ell_1 +1, \ell_4 ^{\prime }}$,     $q\in \overline{\ell_4
^{\prime} +1, \ell}$. We call the operators $\ptl_p ^-$
{\it down-grading operators}, and the operators $\ptl_{p'}^+$ {\it grading
operators}. Set
\begin{eqnarray}
\label{denote2}
\ptl_q =\ptl_q ^-,\hspace{5mm}  \ptl_p =\ptl_p ^-
+\ptl_p ^+ ,\hspace{5mm}
 \ptl_{p^{\prime}} =\ptl_{p^{\prime}} ^+\hspace{5mm}\mbox{and}\hspace{5mm}
\sptl_r=\ptl_{\ell_4+r},
\end{eqnarray}
for $q\in \overline{1, \ell_1}\bigcup \overline{\ell_4 ^{\prime} +1, \ell}$
, $p\in \overline{\ell_1 +1, \ell_3 ^{\prime} }$,  $p^{\prime} \in
\overline{\ell_3 ^{\prime} +1, \ell_4 ^{\prime} }$ and $r\in\ol{1,\ell_5}$.
Then
$\sptl_r,\,r\in\ol{1,\ell_5}$ are the {\it odd} (or
{\it Grassmannian}) {\it derivations\/}
satisfying $\sptl_r^2=0$.

Denote
$\DD=\sum_{p=1} ^\ell \F\ptl_p $. Let $\F[\DD]$ be the (super) polynomial
algebra of $\DD$ with basis
\begin{equation}
\label{1.11}
\mbox{$
\biggl\{\ptl^\mu=\prod\limits_{p=1} ^\ell
\ptl_p ^{\mu_p}\ \biggl|\  \mu=(\mu_1,...,\mu_\ell)\in
\vec K\biggr\},\mbox{ \ \ where \ \ }\vec K=\Z_+ ^{\ell'_4}\times\{0,1\}^{\ell_5},
$}
\end{equation}
and where  $\ptl^{\mu} =1$ if $\mu
=0$. Then the vector space
\begin{eqnarray}
\WW=\WW(\ell_1, \cdots, \ell_5, \G ) =\AA\otimes \F[\DD] ={\rm span}
\{x^{\a, \vec{k}}\ptl^\mu\, |\,(\a ,\vec{k},\mu) \in
\G\times \vec{J}\times\vec K \},
\end{eqnarray}
becomes a $\Z_2$-graded associative superalgebra,
called a {\it superalgebra of generalized Weyl type}, under the operations
(cf.~[SZ1, SZZ])
\begin{eqnarray}
\label{Prod}
u\ptl^\mu \cdot v\ptl^\nu
=u\mbox{$\sum\limits_{\lambda \in \vec K}$} \biggl(
{\sc\,}^{\dis\mu}_{\dis\lambda} \biggr)(-1)^{g(v)g(\mu -\lambda)
+\sum_{\ell'_4< p<q\le\ell}(\mu_q-\l_q)\nu_p}
\ptl^{\lambda }(v)
\ptl^{\mu+\nu -\lambda},
\end{eqnarray}
where $v$ is a homogeneous element of $\AA$ with degree $g(v)\in\Z_2$
(the gradation of the elements in $\WW$ is defined
by (\ref{Gradation-of-A}) and by $g(\ptl_p)=\bar0$ if $p\le\ell'_4$
and $g(\ptl_p)=\bar1$ otherwise), and in general
$g(\mu)=\sum_{p=\ell'_4+1}^\ell\mu_p$ for $\mu\in\vec K$,
and $\bigl({\sc\,}^\mu_\l{\sc\,}\bigr)=
\prod_{p=1} ^\ell \bigl({\sc\,}^{\mu_p}_{\l_p}{\sc\,}\bigr)$
(here $({\sc\,}^j_i{\sc\,})$ is defined to be
$j(j-1)\cdots(j-i+1)/i!$ if $i\ge0$ or $0$ otherwise),
and $\ptl^{\lambda} (v) =
\ptl_1^{\l_1}(\ptl_2^{\l_2}(\cdots(\ptl_\ell^{\l_\ell}(v))))$.
Here we have again used Convention 1.1; for instance, $\ptl^{2_{[\ell'_4+1]}}=0$ since $2_{[\ell'_4+1]}\notin\vec K$.
We shall ALWAYS omit the product notion ``$\cdot$'' when the context is clear.

Under the usual bracket, the superalgebra $\WW$ induces a Lie superalgebra, also denoted by $\WW$ and called a {\it Lie superalgebra of generalized Weyl type}, which is {\it central simple} in the sense that $[\WW,\WW]/\F$ is simple (see [SZZ]).
\vskip10pt
\cl{\bf\S2. \ Preliminaries}\setcounter{section}{2}\setcounter{equation}{0}

Since we shall be mainly interested in infinite dimensional cases, we  assume
$\ell_4 ^{\prime}\ne0$.
Choose a basis of $\WW$ to be
$
B = \{x^{\a, \vec{k}} \ptl^{\mu}\, |\,(\a, \vec{k}, \mu )
\in \G \times \vec{J} \times \vec K\}.
$ 
Fix an element
\begin{eqnarray}
\tau \in\G\mbox{ \ \ such that \ }
\tau_p\neq 0\ \  \mbox{for all }\hspace{2mm} p\in
\overline{\ell_1 +1,\ell_4^{\prime}}.
\end{eqnarray}
For $\mu =(\mu_1, \cdots, \mu_\ell) \in \vec K$, we define the {\it level}
of $\mu$, by $|\mu|=\sum_{p=1} ^\ell \mu_p$, and define a total order on
$\vec K$ by
\begin{eqnarray}
\label{2.3}
\mu < \mu^{\prime} \Leftrightarrow |\mu|<|\mu^{\prime}|,\ {\rm or}\
|\mu| = |\mu^{\prime}|, \,\exists\,p\mbox{ such that }\mu_q=\mu'_q\,(q<p)\mbox{ and }
\mu_p<\mu'_p.
\end{eqnarray}

  Recall that a {\it 2-cocycle} on Lie superalgebra $\WW$ is an $\F$-bilinear
function $\psi: \WW \times \WW\rightarrow \F$, satisfying the following
conditions\begin{eqnarray}
&\!\!\!\!\!\!\!\!\!\!\!\!&\!\!\!\!\!\!\!\!\!\!\!\!
\psi (u_1, u_2) = -(-1)^{g(u_1)g(u_2)} \psi (u_1,u_2)\mbox{ \ (super-skew-symmetry)},\\
\label{2.4}
&\!\!\!\!\!\!\!\!\!\!\!\!&\!\!\!\!\!\!\!\!\!\!\!\!
\psi(u_1,[u_2, u_3]) =\psi([u_1, u_2], u_3) +(-1)^{g(u_1)g(u_2)}\psi(u_2, [u_1 ,u_2])\mbox{ \ (super-Jacobi identity)},
\end{eqnarray}
for $u_1,u_2,u_3 \in \WW$ (we shall always assume an element in $\WW$ is homogeneous).
The vector space of 2-cocycles on $\WW$ is denoted by $C^2 (\WW, \F)$.
A 2-cocycle $\psi$ is called {\it 2-coboundary} or {\it trivial} if there exists an
$\F$-linear function $f$ on $\WW$ such that $\psi=\psi_f$, where
$$\psi_f (u_1, u_2) = f([u_1,u_2])\mbox{ \ \ for \ }u_1,u_2\in\WW.$$
Denote the vector space of all 2-coboundaries by $B^2 (\WW,
\F)$. Two 2-cocycles $\psi, \phi$ are {\it equivalent} if
$\phi -\psi$ is trivial, the quotient space
$$H^2 (\WW ,\F)=C^2
(\WW ,\F)/B^2 (\WW ,\F)$$
 is called the {\it 2-cohomology group} of $\WW$.

The first preliminary result in obtaining our main theorem is the following
rather technical lemma.

{\bf Lemma 2.1.} \ {\it Let $\psi$ be a 2-cocycle on Lie superalgebra
$\WW$,  then there exists a 2-cocycle $\varphi$ equivalent to $\psi$ such
that
\begin{eqnarray}
\label{2.6}
\varphi(t_p\ptl_p, x^{\a, \vec{k}}
\ptl^{\mu})
=0&& \mbox{if \ } \   p\in \overline{1, \ell_3 ^{\prime}},
\\
\label{2.7}
\varphi(\ptl_p, x^{\a, \vec{k}}\ptl^{\mu}) =0
&& \mbox{if \ } \  p\in
\overline{1, \ell},
\\
\label{2.8}
\varphi(s_p\sptl_p, x^{\a, \vec{k}}\ptl^{\mu})
=0
&& \mbox{if \ } \  p\in \overline{1, \ell_5},
\end{eqnarray}
for $(\a, \vec{k }, \mu) \in \G \times \vec{J} \times \vec K$
$($cf.~notations $(\ref{denote1})$ and $(\ref{denote2})).$}

{\it Proof.} \ Define an $\F$-linear function $f$: $\WW\rightarrow \F$ as
follows: For $x^{\a, \vec{k}} \ptl^\mu \in B$ with
$\a \neq 0$, let $q\in\ol{\ell_1+1,\ell'_4}$ be the minimal index such that
$\a_q\neq 0$, we define $f(x^{\a, \vec{k}} \ptl^\mu)$,
inductively on $|k_q|$, by
\begin{eqnarray}
\label{2.9}
 f(x^{\a, \vec{k}} \ptl^\mu) = \left\{{ \begin{array}{ll}
 \a_q ^{-1}(\psi(\ptl_q, x^{\a,\vec{k}}\ptl^\mu)
  - k_q f(x^{\a, \vec{k} -1_[q]}\ptl^\mu) ) & {\rm if} \ k_q\geq 0 ,\\
 -(1+\mu_q)^{-1}(\psi(t_q \ptl_q, x^{\a, \vec{k}}\ptl^\mu)
 -\a_q f(x^{\a, \vec{k}+1_{[q]}}\ptl^\mu ) ) & {\rm if} \ k_q =-1, \\
 (k_q +1)^{-1} (\psi(\ptl_q, x^{\a, \vec{k} +1_{[1_q]}}\ptl^\mu)
 -\a_q f(x^{\a, \vec{k} +1_{[q]}} \ptl^\mu)) & {\rm if} \ k_q \leq
 -2.
\end{array}}
\right.
\end{eqnarray}
Note that
$\ptl_q\in\WW_{\bar0}$ and
$[\ptl_q, x^{\a, \vec{k}} \ptl^{\mu}] =\a_q
x^{\a, \vec{k}} \ptl^{\mu} + k_q x^{\a, \vec{k}
-1_{[q]}} \ptl^{\mu}.$
If $k_q \geq 1$, then
$(\vec{k}-1_{[q]})\in \vec{J},$
so, $x^{\a, \vec{k} -1_{[q]}} \ptl^{\mu} \in B$;
if $k_q\le-1$ then $q\in\ol{\ell'_2+1,\ell'_3}$
and $x^{\a, \vec{k} +1_{[q]}}\ptl^{\mu} \in B,
t_q\ptl_q \in \WW_{\bar{0}},$
$ [t_q\ptl_q, x^{\a,
\vec{k}}
\ptl^{\mu}] =\a_q x^{\a, \vec{k} +1_{[q]}} \ptl^{\mu} + (k_q -\mu_q) x^{\a, \vec{k}}
\ptl^{\mu}$.
Thus the right-hand side of
(\ref{2.9}) makes sense in all cases.

For $t^{\vec{k}}\ptl^{\mu} \in B$ with
$\vec{k} =\vec{i} +\vec{j}$ and $\vec{i} \neq 0$, let $r$ be the
minimal index such that $i_r \neq 0$, we define
\begin{eqnarray}
\label{2.10}
f(t^{\vec{k}} \ptl^\mu) = \left\{\begin{array}{ll} (i_r
-\mu_r)^{-1} \psi (t_r\ptl_r, t^{\vec{k}}
\ptl^\mu) & {\rm if}\hspace{2mm}  i_r \neq \mu_r, \\
(i_r +1)^{-1} \psi (\ptl_r, t^{\vec{k} +1_{[r]}}
\ptl^\mu) & {\rm if}\hspace{2mm}  i_r =\mu_r.
\end{array}
\right.
\end{eqnarray}
Note that in (\ref{2.10})
since $i_r\neq 0,$ we have $r \leq \ell_3 ^{\prime}$ and $t_r\ptl_r \in \WW_{\bar{0}}$, so the right-hand side of (\ref{2.10})
makes sense.

For $s^{\vec{j}}\ptl^\mu \in B$ with $\vec{j} \neq 0$,  let $r'\in\ol{\ell'_4+1,\ell}$
 be the minimal index such
that $ j_{r'} \neq 0$ (i.e., $j_{r'}=1$), we define
\begin{eqnarray}
f(s^{\vec{j}} \ptl^{\mu}) \!=\! \left\{\!\!\begin{array}{ll}
\psi(s_{r'-\ell'_4}\sptl_{r'-\ell'_4}, s^{\vec{j}} \ptl^{\mu})&\!\!
{\rm if} \ \ 1=j_{r'} \neq u_{r'} =0,  \\
\psi(\ptl_1, t_1 s^{\vec{j}} \ptl^{\mu})&\!\!
{\rm if }\ \ \ell_3 ^{\prime} \neq 0, j_{r'}=\mu_{r'}=1,\\
-\tau_{\ell'_4}^{-1}(\mu_{\ell_4 ^{\prime}} +1)^{-1} \biggl(\psi(x^{\tau}, x^{-\tau}
s^{\vec{j}} \ptl^{\mu +1_{[\ell_4 ^{\prime}]}})\\ \!+\!\!
\sum\limits_{
\lambda \in \vec K, \lambda \neq 0, 1_{[\ell_4 ^{\prime}]}} \!\!
\biggl({}^{\dis \mu\! +\! 1_{[\ell_4 ^{\prime}]}}_{\dis\ \ \ \lambda}\biggr)
%
[\tau]^\l
f(s^{\vec{j}} \ptl^{\mu +1_{[\ell_4 ^{\prime}]}
-\lambda})\biggr) &\!\!  \mbox{if }\ell'_3\!=\!0\!\ne\!\ell'_4,\,j_{r'}\!=\!\mu_{r'}\!=\!1,
\end{array}
\right.
\end{eqnarray}
where in general we denote
\begin{equation}
\label{To-use1}
[\a]^\l=\mbox{$\prod\limits_{p=1}^{\ell'_4}$}\a_p^{\l_p}\mbox{ \ \ for \ \ }
\a\in\G,\,\l\in\vec K.
\end{equation}
Note that $[s_{r'-\ell'_4}\sptl_{r'-\ell'_4}, s^{\vec{j}} \ptl^\mu]
\!=\!s^{\vec{j}} \ptl^{\mu}$ if $1\!=\!j_{r'}\! \neq \!\mu_{r'} \!=\!0$, and
$[\ptl_1, t_1s^{\vec{j}} \ptl^{\mu}]=s^{\vec{j}}
\ptl^{\mu} $ if $\ell_3 ^{\prime} \neq 0$, and
\begin{equation}
\label{To-use2}
\mbox{$
[x^{\tau}, x^{-\tau} s^{\vec{j}} \ptl^{\mu +1_{[\ell_4 ^{\prime}]}} ]
\!=\! -\tau_{\ell'_4}(\mu_{\ell_4 ^{\prime}} +1)s^{\vec{j}} \ptl^{\mu}
\!-\!\sum\limits_{\lambda \in \vec K, \lambda \neq 0, 1_{[\ell_4
^{\prime}]}} \!x^{-\tau} s^{\vec{j}}
\biggl({\,}
^{\dis \mu + 1_{[\ell_4 ^{\prime}]}}_{\dis\ \ \ \ \lambda}\,\biggr)
\ptl^{\lambda} (x^{\tau}) \ptl^{\mu +1_{[\ell_4 ^{\prime}]}
-\lambda}
$}
\end{equation}
if $\ell_3 ^{\prime}\!=\!0\! \neq\! \ell_4 ^{\prime}$. We define
$f(s^{\vec{j}}
\ptl^{\mu})$ by induction on $\mu$ with respect to the
order defined in (\ref{2.3}).

Finally for $\ptl^{\mu}\in B$ with $\mu \in \vec K$, if $\ell_3 ^{\prime}\neq
0$, we define
\begin{eqnarray}
f(\ptl^{\mu}) =\psi(\ptl_1, t_1\ptl^{\mu}),
\end{eqnarray}
and
if $\ell_3 ^{\prime}=0\ne\ell_4 ^{\prime}$, we define
(cf.~(\ref{To-use1}) and (\ref{To-use2}))
\begin{eqnarray}
\label{2.13}
\mbox{$
f(\ptl^{\mu})
\!=\!
-(\tau_{\ell'_4}(\mu_{\ell_4 ^{\prime}}\! +\!1))^{-1}
\biggl(\!\psi(x^{\tau}, x^{-\tau} \ptl^{\mu + 1_{[\ell_4 ^{\prime}]}})
\!+\!\!\!\! \sum\limits_{\lambda\in \vec K, \lambda \neq 0, 1_{[\ell_4
^{\prime}]}}\!\!
\bigl({\,}^{\dis \mu\! +\! 1_{[\ell_4 ^{\prime}]}}_{\dis\ \ \ \lambda}\,\bigr)
[\tau]^\l
f(\ptl^{\mu +1_{[\ell_4 ^{\prime}]}-\lambda} )  \biggr)
$}
\end{eqnarray}
by induction on the order of $\mu$.

Now set
$\varphi = \psi - \psi_f$. For $v=x^{\a,\vec k}\ptl^\mu\in B$, we define
\begin{eqnarray}
\label{2.22}
\left\{\begin{array}{ll}
q=q_v= {\rm min}\{q \in \overline{\ell_1 +1, \ell_4 ^{\prime}}\, |\,\a_q \neq 0
\}&\mbox{if \ \ }\a\ne0,
\vs{4pt}\\
r =r_v={\rm min} \{r \in \overline{1, \ell_3 ^{\prime}} \,|\,i_r \neq 0 \}
&\mbox{if \ \ }\vec i\ne0,
\vs{4pt}\\
r'= r'_v=
{\rm min} \{r' \in \overline{\ell_4 ^{\prime} +1, \ell} \,|\,j_{r'} \neq 0 \}
&\mbox{if \ \ }\vec j\ne0,
\end{array}
\right.
\end{eqnarray}
(when there is confusion, we add subscript $v$ to the notation).
Then by (\ref{2.9})--(\ref{2.13}), we
have
\begin{eqnarray}
\label{2.14}
\varphi(\ptl_q, x^{\a ,\vec{k}}\ptl^{\mu}) =0&&
{\rm if } \ \a \neq 0,
\\
\label{2.15}
\varphi(t_q \ptl_q, x^{\a, \vec{k}}\ptl^\mu) = 0&&{\rm if } \ \a \neq0,\,i_q =-1,
\\
\label{2.16}
\varphi(t_r\ptl_r, t^{\vec{i}} \ptl^{\mu}) =0&&{\rm if } \ i_r \neq \mu_r,
\\
\label{2.17}
\varphi(\ptl_r, t^{\vec{k}} \ptl^{\mu}) =0&&{\rm if } \ i_r\geq 2,\, \mu_r =i_r -1,
\ {\rm or} \ r=1,\, \vec{i}=1_{[1]},
\\
\label{2.18}
\varphi(s_{r'-\ell'_4} \sptl_{r'-\ell'_4}, s^{\vec{j}}\ptl^\mu) =0&&  {\rm if } \ 1=j_{r'} \neq\mu_{r'} =0,
\\
\varphi(\ptl_1, t_1 s^{\vec{j}} \ptl^{\mu}) = 0&& {\rm if } \ \ell_3 ^{\prime} \neq0,\,j_{r'} =\mu_{r'} =1,
\\
\label{2.20}
\varphi(x^{\tau}, x^{-\tau} s^{\vec{j}}\ptl^{\mu
+1_{[\ell_4^{\prime}]}} ) =0&&
{\rm if } \ \ell_3 ^{\prime} =0\ne\ell_4 ^{\prime},\,j_{r'} =\mu_{r'} =1,
\\
\label{2.21}
\varphi(x^{\tau}, x^{-\tau} \ptl^{\mu +1_{[\ell_4
^{\prime}]}})=0&&{\rm if } \ \ell'_3=0\ne\ell'_4.
\end{eqnarray}

 Now we prove the lemma in 4 cases.

{\it Case 1:}  \ $\a\neq0$.
Let $q$ be as in (\ref{2.22}), by (\ref{2.14}) and (\ref{2.15}) we obtain
\begin{eqnarray}
\label{2.23}
0&=&\varphi(\ptl_q, \a_p x^{\a, \vec{k} +1_{[p]}}
\ptl^{\mu} +(k_p -\mu_p) x^{\a, \vec{k}} \ptl^{\mu}
)
=\varphi(\ptl_q, [t_p \ptl_p, x^{\a,
\vec{k}}
\ptl^{\mu}])\nonumber\\
&=&\varphi ([\ptl_q, t_p \ptl_p], x^{\a,
\vec{k}} \ptl^{\mu}) +\varphi (t_p \ptl_p,
[\ptl_q, x^{\a, \vec{k}} \ptl^{\mu}])\nonumber \\
&=&\delta_{p, q}\varphi(\ptl_p, x^{\a, \vec{k}}
\ptl^{\mu}) + \varphi (t_p \ptl_p, \a_q x^{\a, \vec{k}} \ptl^{\mu} +
i_q x^{\a, \vec{k} -1_{[q]}} \ptl^{\mu})\nonumber\\
&=& \a_q \varphi (t_p\ptl_p, x^{\a,
\vec{k}}\ptl^{\mu}) +
 i_q \varphi (t_p\ptl_p, x^{\a, \vec{k} -1_{[q]}}
 \ptl^{\mu})
\mbox{ \ \ for $p\in \overline{1, \ell_3 ^{\prime}}$}.
\end{eqnarray}
We also have
\begin{eqnarray}
\label{2.24}
0&=& \varphi(\ptl_q, \a_p x^{\a, \vec{k}}
\ptl^{\mu} \pm i_p x^{\a, \vec{k} -1_{[p]}}
\ptl^{\mu})
=\varphi(\ptl_q, [\ptl_p, x^{\a, \vec{k}}
\ptl^{\mu}])\nonumber \\
&=&\a_q \varphi (\ptl_p, x^{\a, \vec{k}}
\ptl^{\mu}) +i_q \varphi (\ptl_p, x^{\a, \vec{k}
-1_{[q]}}\ptl^{\mu})
\mbox{ \ \ for $p \in \overline{1, \ell} $,}
\end{eqnarray}
(note that when $p\in\ol{\ell'_4+1,\ell}$, $\ptl_p$ is an odd derivation, so it may produce a minus sign when applying it to $x^{\a,\vec k}$).
If $i_q \geq 0$, then (\ref{2.6}) and
(\ref{2.7}) follow from (\ref{2.23}), (\ref{2.24}) and induction on $i_q$.
If $i_q =-1,\,p=q$, then (\ref{2.6}) and (\ref{2.7}) follow from
(\ref{2.14}) and (\ref{2.15}).
Assume $i_q=-1,\,p\neq q$. Using (\ref{2.15}) we
have
\begin{eqnarray}
\label{2.25}
0&\!\!\!=\!\!\!& \varphi (t_q \ptl_q, [t_p \ptl_p,
x^{\a ,\vec{k}}
\ptl^{\mu}])
=
\varphi(t_p\ptl_p ,[t_q\ptl_q,
x^{\a, \vec{k}} \ptl^{\mu}]) \nonumber \\
&\!\!\!=\!\!\!& \varphi(t_p\ptl_p, \a_q x^{\a, \vec{k}
+1_{[q]}} \ptl^{\mu} +i_qx^{\a, \vec{k}} \ptl^{\mu}
-\mu_q x^{\a, \vec{k}}\ptl^{\mu})
=
-(1 +\mu_q) \varphi(t_p \ptl_p, x^{\a,
\vec{k}}\ptl^{\mu}),
\end{eqnarray}
which implies (\ref{2.6}), where the last equality of (\ref{2.25}) follows from
the fact that (\ref{2.6}) holds when $i_q\ge0$.
Using (\ref{2.23}), the proof of (\ref{2.6}) can be completed by induction on $-i_q$ when $i_q\le-2$.
To complete the proof of (\ref{2.7}), first suppose $p
\in \overline{1, \ell_4 ^{\prime}}$. If $i_q =-1$, then
\begin{eqnarray}
\label{Add1}
0\!\!\!\!&=\!\!\!\!& \varphi (t_q \ptl_q, [\ptl_p, x^{\a, \vec{k}} \ptl^{\mu}])\nonumber\\
&=\!\!\!\!&\varphi([t_q \ptl_q, \ptl_p], x^{\a,
\vec{k}}
\ptl^{\mu}) + \varphi(\ptl_p, [t_q\ptl_q, x^{\a, \vec{k}}\ptl^\mu])
\nonumber
\\
&=\!\!\!\!&-\delta_{p,q} \varphi(\ptl_p, x^{\a, \vec{k}}
\ptl^{\mu})\! +\!\a_q \varphi(\ptl_p, x^{\a, \vec{k} +1_{[q]}} \ptl^{\mu})
-(1\! +\!\mu_q) \varphi(\ptl_p, x^{\a, \vec{k}}\ptl^{\mu})
\nonumber
\\
&=\!\!\!\!&-(1 +\delta_{p, q} +\mu_q) \varphi(\ptl_p,
x^{\a,\vec{k}}\ptl^{\mu}),
\end{eqnarray}
which implies (\ref{2.7}), where the last equality of (\ref{Add1}) follows from
the fact that (\ref{2.7}) holds when $i_q\ge0$. If $i_q\le-2$, (\ref{2.7}) is obtained from
(\ref{2.24}) by induction on $-i_q$.
Next suppose $p \in \overline{\ell_4 ^{\prime}
+1, \ell}$.
If $i_q =-1$ then from (\ref{2.6}) we have (noting that $p\ne q$ in this case)
\begin{eqnarray}
\label{iq=-1}
0\!\!\!\!&=\!\!\!\!& \varphi(t_q\ptl_q, [\ptl_p, x^{\a,
\vec{k}} \ptl^{\mu}])
\nonumber
=\varphi([t_q\ptl_q, \ptl_p], x^{\a,
\vec{k}}
\ptl^{\mu}) +\varphi(\ptl_p, [t_q \ptl_q, x^{\a, \vec{k}} \ptl^{\mu}]) \\
&=\!\!\!\!& \a_q \varphi(\ptl_p, x^{\a, \vec{k} +1_{[q]}}
\ptl^{\mu}) -(1 + \mu_q) \varphi(\ptl_p, x^{\a, \vec{k}} \ptl^{\mu})
=-(1+ \mu_q) \varphi(\ptl_p, x^{\a, \vec{k}}
\ptl^{\mu}),
\end{eqnarray}
which implies (\ref{2.7}).
The proof of (\ref{2.7}) is completed by (\ref{2.24}) and induction on $-i_q$ when $i_q\le-2$.  The proof of (\ref{2.8}) is similar to that of (\ref{2.6}).

{\it Case 2:} \ $\a=0,\,\vec{i} \neq 0$ (which implies $\ell_3 ^{\prime} \geq 1$). Let $r$ be as in
(\ref{2.22}). First we prove (\ref{2.6}). If $i_r \neq\mu_r$, then by
(\ref{2.16}) we have
\begin{eqnarray}
\label{2.29}
0&\!\!\!=\!\!\!&\varphi(t_r \ptl_r, [t_p\ptl_p,
t^{\vec{k}}\ptl^\mu])
=
\varphi([t_r \ptl_r, t_p\ptl_p], t^{\vec{k}} \ptl^{\mu}) +\varphi(t_p\ptl_p, [t_r\ptl_r, t^{\vec{k}} \ptl^{\mu}])\nonumber
\\
&\!\!\!=\!\!\!& (i_r-\mu_r) \varphi(t_p \ptl_p, t^{\vec{k}}
\ptl^{\mu})
\hspace {5mm}
(\mbox{noting that }r, p \in \overline{1,
\ell_3^{\prime}}),
\end{eqnarray}
which implies (\ref{2.6}). Assume $i_r = \mu_r$. Then we have $i_r \geq 1$ (cf.~(\ref{2.22})), and
\begin{eqnarray}
\varphi(t_p \ptl_p, t^{\vec{k}}
\ptl^{\mu}) &\!\!\!=\!\!\!&(1 + i_r)^{-1} \varphi(t_p \ptl_p,
[\ptl_r, t^{\vec{k} + 1_{[r]}}
\ptl^{\mu}])\nonumber\\
&\!\!\!=\!\!\!&(1 +i_r)^{-1} (-\delta_{r,p} +i_p +\delta_{p, r} -\mu_p) \varphi(\ptl_r, t^{\vec{k} +1_{[r]}} \ptl^{\mu})
= 
0,
\end{eqnarray}
which gives (\ref{2.6}),
where the last equality follows from (\ref{2.17}) by noting that the condition of
(\ref{2.17}) is satisfied by $t^{\vec k+1_{[r]}}\ptl^\mu$ since $i_r+1=\mu_r+1\ge2$.
 Similarly we have (\ref{2.8}). Next consider (\ref{2.7}). First suppose $p\! \in\! \overline{1, \ell_4
^{\prime}}$. Since
\begin{equation}
t^{\vec{k}} \ptl^{\mu} =\left\{\begin{array}{ll} (i_1 +1)^{-1}
[\ptl_1 ,t^{\vec{k} +1_{[1]}} \ptl^{\mu}]
&  {\rm if}  \ i_1 \neq-1,\\
-(1 +\mu_1)^{-1} [t_1 \ptl_1,t^{\vec{k}}
\ptl^{\mu}]  & {\rm if} \ i_1 =-1,
\end{array}
\right.
\end{equation}
we have
\begin{equation}
\varphi(\ptl_p ,t^{\vec{k}} \ptl^{\mu})
=\left\{\begin{array}{ll} (i_1 +1)^{-1}
(i_p+\delta_{p,1})\varphi(\ptl_1 ,t^{\vec{k} +1_{[1]}
-1_{[p]}}
\ptl^{\mu}) & {\rm if} \ i_1 \neq -1, \\
-(1+ \mu_1)^{-1}\d_{p,1} \varphi(\ptl_1 ,t^{\vec{k}} \ptl^{\mu}) &
{\rm if} \ i_1 =-1, \end{array}
\right.
\end{equation}
where the second case
is obtained from the super-Jacobi identity and (\ref{2.6}).
So the proof of (\ref{2.7}) is reduced the case $p=1$. Using (\ref{2.6}) we have
\begin{equation}
 0= \varphi(t_1 \ptl_1, [\ptl_1, t^{\vec{k}}
\ptl^{\mu}]) =(-1 +i_1 -\mu_1) \varphi(\ptl_1, t^{\vec{k}}
\ptl^{\mu}).
\end{equation}
Thus it remains to consider the case $i_1 =\mu_1 +1$. The
result follows from (\ref{2.17}) if $\mu_1 \geq 1$ or $\vec i =1_{[1]}$. Thus assume
$\mu_1 =0$ and $\vec i \neq 1_{[1]}$. Then $i_{q'} \neq 0$ for some $q'
\in \overline{2, \ell_3  ^{\prime}}$ and we denote $q'$ to be the minimal index with $i_{q'} \neq  0$. If $i_{q'} \neq \mu_{q'}$, then by (\ref{2.6}) we have
\begin{equation}
0 =\varphi(t_{q'} \ptl_{q'} ,[\ptl_1 ,t^{\vec{k}}
\ptl^{\mu}]) =(i_{q'} -\mu_{q'}) \varphi(\ptl_1 ,t^{\vec{k}}
\ptl^{\mu}),
\end{equation}
and if $i_{q'} =\mu_{q'}$ then
\begin{equation}
\varphi(\ptl_1 ,t^{\vec{k}} \ptl^{\mu}) =(i_{q'} +1)^{-1}
\varphi(\ptl_1 ,[\ptl_{q'} ,t^{\vec{k} +1_{[q']}}
\ptl^{\mu}])
=(i_{q'} +1)^{-1} \varphi(\ptl_{q'} ,t^{\vec{k}
-1_{[1]} +1_{[q']}}
\ptl^{\mu}) =0,
\end{equation}
where the last equality follows from (\ref{2.17}) by noting that $q'$ is precisely
the number
$r_v$ defined in (\ref{2.22}) for $v=t^{\vec{k}-1_{[1]} +1_{[q']}}\ptl^{\mu}$.
Finally assume $p \in \overline{\ell_4 ^{\prime} +1,\ell}$.
From the result above we
have
\begin{equation}
0 = \varphi(\ptl_1 ,[\ptl_p ,t^{{\vec{k}} +1_{[1]}}
\ptl^{\mu}]) =\varphi(\ptl_p ,[\ptl_1 ,t^{\vec{k} +1_{[1]}} \ptl^{\mu}]) =(i_1 +1) \varphi(\ptl_p ,t^{\vec{k}}
\ptl^{\mu}),
\end{equation}
which implies (\ref{2.7}) if $i_1 \neq-1$.
If $i_1 =-1$, then (\ref{2.6}) gives
\begin{equation}
0= \varphi(t_1 \ptl_1 ,[\ptl_p ,t^{\vec{k}}
\ptl^{\mu} ]) =\varphi(\ptl_p ,[t_1 \ptl_1 ,t^{\vec{k}} \ptl^{\mu}]) =-(1 + \mu_1) \varphi(\ptl_p ,t^{\vec{k}}
\ptl^{\mu}),
\end{equation}
which implies (\ref{2.7}) since $\mu_1\ge0$.

{\it Case 3:} \ $\a =0,\, \vec{i} =0, \,\vec{j} \neq 0.$
Consider (\ref{2.6}).
For $p \in \overline{1, \ell_3 ^{\prime} }$, we have
\begin{equation}
\label{Case3.1}
\varphi(t_p \ptl_p ,s^{\vec{j}} \ptl^{\mu})
=\varphi(t_p\ptl_p,[\ptl_p,t_p s^{\vec{j}}\ptl^{\mu}]) =0,
\end{equation}
where the last equality follows from the super-Jacobi identity and Case 2.
Consider (\ref{2.8}). Assume $p\in\ol{1,\ell_5}$. If $\ell_3 ^{\prime} \neq 0$,
then similar to (\ref{Case3.1}), we have
$$\varphi(s_p \sptl_p ,s^{\vec{j}} \ptl^{\mu})
=\varphi(s_p\sptl_p,[\ptl_1,t_1s^{\vec{j}}\ptl^{\mu}]) =0.$$
Assume $\ell_3 ^{\prime} =0\ne\ell_4 ^{\prime}$.
Let $r'$ be as in (\ref{2.22}). If $\mu_{r'}=0$, we have
$$\varphi(s_p\sptl_p,s^{\vec j}\ptl^\mu)
=\varphi(s_p\sptl_p,[s_{r'-\ell'_4}\sptl_{r'-\ell'_4},s^{\vec j}\ptl^\mu])=0
\mbox{ \ (cf.~notations (\ref{denote1}) and (\ref{denote2}))},
$$
where the last equality follows from super-Jacobi identity and (\ref{2.18}).
If $\mu_{r'}\!=\!1\!=\!j_{r'}$, we have
\begin{eqnarray}
0& =&\varphi([s_p\sptl_s ,x^{\tau}] ,x^{-\tau}
s^{\vec{j}}
\ptl^{\mu + 1_{[\ell_4 ^{\prime}]}})\mbox{\ \ \ \
(since $[s_p\sptl_p , x^{\tau}] =0$)}\nonumber\\
&=& \varphi(s_p \sptl_p ,[x^{\tau} ,x^{-\tau}
s^{\vec{j}}
\ptl^{\mu +1_{[\ell_4 ^{\prime}]}}])\mbox{\ \ \ \
(from (\ref{2.4}) and (\ref{2.20}))}\nonumber\\
&=&-\mbox{$\sum \limits_{0 \neq \lambda \in \vec K} $}
\biggl(\,^{\dis \mu +1_{[\ell_4 ^{\prime}]}}_{\dis
\ \ \ \lambda}\biggr)
[\tau]^\l
\varphi(s_p \sptl_p , s^{\vec{j}} \ptl^{\mu
-\lambda +1_{[\ell_4 ^{\prime}]}})
\mbox{\ \ \ \ (cf.~(\ref{To-use1}))},
\end{eqnarray}
and (\ref{2.8}) follows from the induction on $|\mu|$.
Consider (\ref{2.7}). If $p \in \overline{1, \ell_3 ^{\prime}}$,
we have
\begin{eqnarray}
\label{2.38}
0&=&\varphi(t_p \ptl_p ,[\ptl_p ,s^{\vec{j}}\ptl^{\mu}])
\hspace{0.5cm} \mbox{\rm (since \ $[\ptl_p ,s^{\vec{j}} \ptl^{\mu}] =0)$}
\nonumber
\\
&=& \varphi([t_p \ptl_p ,\ptl_p] ,s^{\vec{j}}
\ptl^{\mu}) +\varphi(\ptl_p ,[t_p \ptl_p ,s^{\vec{j}} \ptl^{\mu}])
=
-(1+ \mu_p) \varphi(\ptl_p ,s^{\vec{j}}
 \ptl^{\mu}).
\end{eqnarray}
Assume $p \in \overline{\ell_3 ^{\prime}+1, \ell}$. If $\ell_3 ^{\prime} \geq
  1$, we have
\begin{equation}
\varphi(\ptl_p, s^{\vec{j}}
\ptl^{\mu})=\varphi(\ptl_p, [\ptl_1, t_1
s^{\vec{j}}
\ptl^{\mu}]) =\varphi(\ptl_1, [\ptl_p, t_1 s^{\vec{j}}
\ptl^{\mu}]) =0\hspace{0.4cm}\mbox{ (by Case 2)}.
\end{equation}
Assume $\ell'_3=0$.
If there exists $p'\in\ol{\ell'_4+1,\ell}{\,}\bs\{p\}$ such that
$j_{p'}\ne\mu_{p'}$, then
$$
\varphi(\ptl_p,s^{\vec j}\ptl^\mu)
=\varphi(\ptl_p,[s_{p'-\ell'_4}\sptl_{p'-\ell'_4},s^{\vec j}\ptl^\mu])=0
\mbox{ \ (by (\ref{2.4}) and (\ref{2.8}))}.
$$
Assume
\begin{equation}
\label{Assume1}
\mbox{
$j_{p'}=\mu_{p'}$ for all $p'\in\ol{\ell'_4+1,\ell}{\,}\bs\{p\}$.
}
\end{equation}
If $p=r'$ and $1=j_{r'}=\mu_{r'}$, we have
$$
0=\varphi(s_{r'-\ell'_4}\sptl_{r'-\ell'_4},[\ptl_p,s^{\vec j}\ptl^\mu])
=-\varphi(\ptl_p,s^{\vec j}\ptl^\mu)\mbox{ \ (cf.~notations (\ref{denote1})
and (\ref{denote2}))}
$$
by super-Jacobi identity. Assume $p=r'$ and $1=j_{r'}\ne\mu_{r'}=0$. Let
$v=[\ptl_p ,x^{-\tau} s^{\vec{j}} \ptl^{\mu +1_{[\ell_4
^{\prime}]}}]=x^{-\tau} s^{\vec{j}-1_{[r']}} \ptl^{\mu +1_{[\ell_4
^{\prime}]}}$. Then $v$ either has the form $x^{-\tau}\ptl^{\mu +1_{[\ell_4
^{\prime}]}}$ (if $s^{\vec j}=1_{[r']}$) or $r'_v\ne r'$ satisfies $j_{r'_v}=\mu_{r'_v}=1$ by (\ref{Assume1}) and by definition (\ref{2.22}); in either
case, we have
\begin{equation}
\label{more0}
\varphi(x^{\tau},v)=0\mbox{ \ (by (\ref{2.20}) or (\ref{2.21}))}.
\end{equation}
Thus
\begin{eqnarray}
\label{2.40}
0 \!\!\!&=&\!\!\!\varphi([\ptl_p ,x^{\tau}] ,x^{-\tau} s^{\vec{j}}
\ptl^{\mu +1_{[\ell_4 ^{\prime}]}} )
\mbox{\ \ \ \ (since $[\ptl_p ,x^{\tau}]=0$)}
\nonumber\\& =&\!\!\!
\varphi(\ptl_p ,[x^{\tau} ,x^{-\tau} s^{\vec{j}}
\ptl^{\mu +1_{[\ell_4 ^{\prime}] }} ] ) +\varphi([\ptl_p ,x^{-\tau} s^{\vec{j}} \ptl^{\mu +1_{[\ell_4
^{\prime}]}}],
x^{\tau})\mbox{\ \ \ \ (by (\ref{2.4}))}
\nonumber\\
&=&\!\!\!\varphi(\ptl_p ,[x^{\tau} ,x^{-\tau} s^{\vec{j}}
\ptl^{\mu +1_{[\ell_4 ^{\prime}]} }])
\nonumber\\
 &=&\!\!\!-\mbox{$\sum \limits_{0 \neq\lambda \in \vec K}$}
\biggl(\,^{\dis\mu +1_{[\ell_4^{\prime}]}}
_{\dis\ \ \ \lambda} \biggr)
[\tau]^\l
\varphi(\ptl_p ,s^{\vec{j}}\ptl^{\mu +1_{[\ell_4
^{\prime}]} -\lambda} )\mbox{\ \ \ \ (cf.~(\ref{To-use2}))},
\end{eqnarray}
which gives (\ref{2.7}) by induction on $\mu$, where the third equality follows from
(\ref{more0}). Finally assume $p\ne r'$.
Then by (\ref{Assume1}) we must have $1=j_{r'}=\mu_{r'}$, and we still have (\ref{2.40}) where now the third equality follows from (\ref{2.20}).

{\it Case 4:} \ $\a =0,\, \vec{k} =0.$
From (\ref{Case3.1}), we have (\ref{2.6}).
Consider (\ref{2.7}). If $p \in \overline{1, \ell_3 ^{\prime}},$ then
\begin{equation}
0 = \varphi(t^{1_{[p]}} \ptl_p, [\ptl_p, \ptl^{\mu}])
=-(1 +\mu_p)\varphi(\ptl_p, \ptl^{\mu}),
\end{equation}
which implies (\ref{2.7}).
Assume $p \in\overline{\ell_3 ^{\prime} +1, \ell}$.
If $\ell_3 ^{\prime} \geq 1$, then
\begin{equation}
\varphi(\ptl_p, \ptl^{\mu}) =\varphi(\ptl_p,
[\ptl_1, t^{1_{[1]}} \ptl^{\mu} ]) =\varphi(\ptl_1,
[\ptl_p, t^{1_{[1]}} \ptl^{\mu}]) =0.
\end{equation}
If $\ell_3 ^{\prime} =0$, then
\begin{equation}
\begin{array}{lll}
0 &=& \varphi([\ptl_p ,x^{\tau } ], x^{-\tau} \ptl^{\mu +
1_{\ell_4 ^{\prime}}})
=
- \sum \limits_{0 \neq \lambda  \in \vec K} \biggl(\,^{\dis\mu +1_{\ell_4 ^{\prime}}}_{\dis\ \ \ \lambda}\biggr) \prod\limits_q \tau_q ^{\lambda_q} \varphi(\ptl_p, \ptl^{\mu +1_{[\ell_4 ^{\prime}]} -\lambda
})
\end{array}
\end{equation}
by (\ref{To-use2}), (\ref{2.21}) and by noting that
if $p \leq \ell_4 ^{\prime}$ then $[\ptl_p, x^{\tau}]
=\tau_p x^{\tau}, \,[\ptl_p, x^{-\tau} \ptl^{\mu +1_{[\ell_4 ^{\prime}]}}]
  =-\tau_p x^{-\tau} \ptl^{\mu +1_{[\ell_4 ^{\prime}]}}$, and if $p > \ell_4^{\prime}$ then $ [\ptl_p, x^{\tau}] =0,\, [\ptl_p, x^{-\tau} \ptl^{\mu + 1_{[\ell_4
 ^{\prime}]}}] =0$. Induction on $|\mu|$ gives (\ref{2.7}).
Finally consider (\ref{2.8}). Suppose $p\in\ol{1,\ell_5}$.
If $\ell_3 ^{\prime} \neq 0$, we have
\begin{equation}
0=\varphi(\ptl_1 ,[s_p\sptl_p ,t_1\ptl^{\mu} ])
=\varphi(s_p\sptl_p,[\ptl_1,t_1\ptl^{\mu}])
=\varphi(s_p\sptl_p,\ptl^{\mu}).
\end{equation}
If $\ell_3 ^{\prime} =0 ,\, \ell_4 ^{\prime} \neq 0$ then by (\ref{2.21}) and (\ref{To-use2}) we have
$$ 
\begin{array}{lll}
0&\!\!\!=\!\!\!&\varphi([s_p\sptl_p , x^{\tau}] ,x^{-\tau}
\ptl^{\mu +1_{[\ell_4 ^{\prime}]}})
=
\varphi(s_p\sptl_p , [x^{\tau} ,x^{-\tau}
\ptl^{\mu +1_{[\ell_4 ^{\prime}]}}]) + \varphi(x^{\tau}
,[s_p \sptl_p ,x^{-\tau}
\ptl^{\mu +1_{[\ell_4 ^{\prime}]}}]) \\
&\!\!\!=\!\!\!&\sum \limits_{0 \neq \lambda \in \vec K}
 \biggl(\,^{\dis\mu +1_{[\ell_4 ^{\prime}]}}_{\dis\ \ \ \lambda}\biggr)
[\tau]^\l\varphi(s_p\sptl_p ,\ptl^{\mu -\lambda +1_{[\ell_4 ^{\prime}]}}
),
\end{array}
$$ 
which gives the result by
induction on $|\mu|$. This completes the proof of the lemma.
\hfill$\Box$\vskip10pt
\cl{\bf \S3. \ Main \vs{-0.1pt}results}
\setcounter{section}{3}
\setcounter{equation}{0}
Recall that we assume $\ell'_4>0$.
Denote
\begin{eqnarray}
\label{3.02}
\WW_0\!=
\!{\rm span}\{x^{\a,\vec i}\ptl^\mu\,|\,(\a,\vec i,\mu)
\!\in\!\G\!\times\!\vec J'\!\times\!
\vec K'\}, \
\vec J'\!=\!\Z_+^{\ell'_2}\!\times\!\Z^{\ell_3}\!\times\!\{0\}^{\ell_4+\ell_5},\
\vec K'\!=\!\Z_+^{\ell'_4}\!\times\!\{0\}^{\ell_5}
\end{eqnarray}
(cf.~(\ref{1.3}) and (\ref{1.11})).
Then $\WW_0$ is a Lie algebra of Weyl type whose $2$-cocycles were
considered in [S2] (cf.~also [L], [LW], [S1]).
The following result can be found in [S2].

{\bf Theorem 3.1.} (1) {\it If $\ell'_4=\ell_4=1$, then
$H^2(\WW_0,\F)=\F\ol\phi_0$, where $\ol\phi_0$ is the
cohomology class of $\phi_0$ defined by
\begin{equation}
\label{Lem3.1}
\phi_0(x^\a[\ptl_1]_\mu,x^\b[\ptl_1]_\nu)
=\d_{\a+\b,0}(-1)^\mu\mu!\nu!(^{\ \dis\ \a+\mu}_{\dis \mu+\nu+1}),
\end{equation}
for $\a,\b\in\G\subseteq\F,\,\mu,\nu\in\Z_+,$ where $[\ptl_1]_\mu=\ptl_1(\ptl_1-1)\cdots(\ptl_1-\mu+1)$.}
\par
(2) {\it If $\ell'_4=\ell_3=1$, then for any $\g\in\G$, there exists a cohomology class
$\ol\phi_\g\in H^2(\WW_0,\F)$ defined by
\begin{equation}\label{Lem3.1.2}
\phi_\g(x^{\a,i}\ptl^\mu,x^{\b,j}\ptl^\nu)
=\d_{\a+\b,\g}(-1)^\mu\mu!\nu!
\mbox{$\sum\limits_{s=0}^{\mu+\nu+1}$}
\bigl({\sc\,}^{\ssc\,\dis i}_{\dis s}{\sc\,}\bigr)
\frac{\a^{\mu+\nu+1-s}}{(\mu+\nu+1-s)!}
\cdot\frac{\g^{s-i-j-1}}{(s-i-j-1)!}\ ,
\end{equation}
for all $(\a,i,\mu),(\b,j,\nu)\in\G\times\Z\times\Z_+$,
where as in [S2], $\frac{1}{k!}$ is understood as zero when $k<0$, and
when $\a$ is taken value $0$, it is understood as $\lim_{\a\to0}\a$ (thus in particular, $\a^{\mu +\nu +1 -r}=1$ if ${\mu +\nu +1 -r}=0$ and $\a=0$).
Furthermore $H^2(\WW_0,\F)=\prod_{\g\in\G}\F\ol\phi_\g$ is a direct product.
\par
{\rm(3)} If $\ell_1+\ell_2\ge1$ or $\ell'_4\ge2$, then $H^2(\WW_0,\F)=0.$
\hfill$\Box$
}

{\bf Remark 3.2.} \
(1) It is proved in [S2] that $\phi_\g$ is in fact a $2$-cocycle of the associative
algebra $\WW_0$, satisfying
\begin{equation}
\label{Rmk3.2}
\phi_\g(a,bc)+\phi_\g(b,ca)+\phi_\g(c,ab)=0\mbox{ \ for \ }a,b,c\in\WW_0.
\end{equation}

\def\dx{\mbox{$\frac{d}{dx}$}}
(2) In the case of Theorem 3.1(1), we can suppose $1\in\G$ (see e.g., [SZ2]).
Then we can define the {\it derivative} $\dx$ by $\dx x^\a=\a x^{\a-1}$. We have
$[\ptl_1]_\mu=x^\mu(\dx)^\mu$ and (\ref{Lem3.1.2}) becomes
\begin{equation}
\label{Rmk3.2.2}
\phi_0(x^{\a+\mu}(\dx)^\mu,x^{\b+\nu}(\dx)^\nu)=
\d_{\a+\b,0}(-1)^\mu\mu!\nu!(^{\ \dis\ \a+\mu}_{\dis \mu+\nu+1})\mbox{ \ for \ }
\a,\b\in\G,\,\mu,\nu\in\Z_+
\end{equation}
(this $2$-cocycle for the case $\G=\Z$ (the classical Weyl algebra) seems to appear first in [KP]).
We prove as follows that $\phi_0$ also satisfies (\ref{Rmk3.2}):
First by (\ref{Rmk3.2.2}), we have
\begin{equation}
\label{Rmk3.2.21}
\phi_0(x^{\a+\mu}(\dx)^\mu,x^{\b+\nu}(\dx)^\nu)=
(-1)^{\nu+1}\frac{\mu!\nu!}{(\mu+\nu+1)!}\phi_0(x^1,x^{\a+\mu}(\dx)^{\mu+\nu+1}(x^{\b+\nu})).
\end{equation}
Using this and $x^{\a+\mu}(\dx)^\mu x^{\b+\nu}(\dx)^\nu=\sum_{\l\in\Z_+}(^\mu_\l)
[\b+\nu]_\l x^{\a+\b+\mu+\nu-\l}(\dx)^{\mu+\nu-\l}$ (where
$[\b+\nu]_\l$ is a similar notation to $[\ptl_1]_\l$, cf.~(\ref{Lem3.1}), (\ref{Prod})), we have
\begin{equation}\!
\label{Rmk3.2.22}
\phi_0(a(\dx)^\mu,b(\dx)^\nu c(\dx)^\l)\!
=\!\!\mbox{$\sum\limits_{s=0}^\nu$}(-1)^{\nu+\l+1-s}
\frac{\mu!(\nu\!+\!\l\!-\!s)!}{(k\!+\!1\!-\!s)!}
(^{\dis\nu}_{\ssc\,\dis s})
\phi_0(x^1,a(\dx)^{k+1-s}(b(\dx)^s(c))),\!\!\!
\end{equation}
for $a,b,c\in\AA_{\bar0}$
(cf.~(\ref{Gradation-of-A})), where $k=\mu+\nu+\l$. Using shifted version of
(\ref{Rmk3.2.22}), we have
\begin{eqnarray}
\label{Rmk3.2.23}
\!\!\!\!\!\!\!\!\!&\!\!\!\!\!\!\!\!\!&\!\!\!\!\!\!\!\!\!
\phi_0(b(\dx)^\nu,c(\dx)^\l a(\dx)^\mu)
=
\mbox{$\sum\limits_{s=0}^\l$}
(-1)^{\nu+s}
\frac{\nu!(\l\!+\!\mu\!-\!s)!}{(k\!+\!1\!-\!s)!}
(^{\dis\l}_{\ssc\,\dis s})
\phi_0(x^1,a(\dx)^s((\dx)^{k+1-s}(b)c)),
\\
\label{Rmk3.2.24}
\!\!\!\!\!\!\!\!\!&\!\!\!\!\!\!\!\!\!&\!\!\!\!\!\!\!\!\!
\phi_0(c(\dx)^\l,a(\dx)^\mu b(\dx)^\nu)
=
\mbox{$\sum\limits_{s=0}^\mu$}
(-1)^\l
\frac{\l!(\nu\!+\!\mu\!-\!s)!}{(k\!+\!1\!-\!s)!}(^{\dis\mu}_{\ssc\,\dis s})
\phi_0(x^1,a(\dx)^s(b)(\dx)^{k+1-s}(c)).
\end{eqnarray}
Denote the right-hand sides of (\ref{Rmk3.2.22})--(\ref{Rmk3.2.24}) by
$
\sum_{s=0}^{k+1}d_{p,s}\phi(t^1,a(\dx)^s(b)(\dx)^{k+1-s}(c))
\mbox{ for }p=1,2\mbox{ and 3 respectively.}
$
Using
$(1+x)^{k+1-s}(1+x)^{-(\l+1)}=(1+x)^{\mu+\nu-s}$, we deduce the
binomial formula
$
\sum_q(-1)^q(^{\dis k+1-s}_{\dis\ \ \nu-q})(^{\dis\l+q}_{\,\dis\ \ q})=
(^{\dis\mu+\nu-s}_{\,\dis\ \ \ \ \,\nu}).
$
From this, we can deduce that if $s\le\mu$,
then
$$
\begin{array}{lll}
d_{1,s}\!\!\!\!&=\sum\limits_{q=0}^\nu\dis
(-1)^{\nu+\l+1-q}
{\mu!(\nu\!+\!\l\!-\!q)!\over(k+1-q)!}
(^{\dis\nu}_{\ssc\,\dis q})
(^{\dis k\!+\!1\!-\!q}_{\dis\!\!\ \ \ \ \ s})
=(-1)^{\l+1}{\l!(\nu\!+\!\mu\!-\!s)!\over(k+1-s)!}(^{\dis\mu}_{\ssc\,\dis s})
=-d_{3,s},
\hfill\cr
\end{array}
$$
and $d_{2,s}=0$; and if $\mu<s\le\mu+\nu$, then $d_{1,s}=d_{2,s}=d_{3,s}=0$;
and if $\mu+\nu<s\le k+1$,
then
$$
\begin{array}{lll}
d_{1,s}\!\!\!\!&
=\sum\limits_{q=0}^\nu\dis
(-1)^{\nu+\l+1-q}
{\mu!(\nu+\l-q)!\over(k+1-q)!}
(^{\dis\nu}_{\ssc\,\dis q})
(^{\dis k+1-q}_{\dis\ \ \ \ \ s})
\vs{4pt}\hfill\cr&
=
\sum\limits_{q=0}^\l\dis(-1)^{\nu+q+1}
{\nu!(\l+\mu-q)!\over(k+1-q)!}
(^{\dis\l}_{\ssc\,\dis q})
(^{\dis\ \ \ \ \ q}_{\dis k+1-s})
=-d_{2,s},
\end{array}
$$
and $d_{3,s}=0$.
This proves that the sum of (\ref{Rmk3.2.22})--(\ref{Rmk3.2.24}) is zero.
\hfill$\Box$

Assume $\ell_5>0$. Observe that as an associative superalgebra under the product
(\ref{Prod}), $\WW$ can be decomposed into the following tensor product of
super-subalgebras:
\begin{eqnarray}
\label{3.01}
\WW=\WW_0\otimes\WW_1,\ \ \
\WW_1=\WW_{1,1}\otimes\WW_{1,2}\otimes\cdots\otimes\WW_{1,\ell_5},
\end{eqnarray}
where
\begin{eqnarray}
\label{3.03}
\WW_{1,p}\!=\!{\rm span}\{s_p^m\sptl_p^n\,|\,m,n\!\in\!\{0,1\}\}
\mbox{ \ (a superalgebra of
 dimension $4$) \ for $p\!\in\!\ol{1,\ell_5}$}
\end{eqnarray}
(cf.~notations (\ref{denote1}) and (\ref{denote2})).
Denote
\begin{eqnarray}
\label{3.04}
\WW^\wedge_p=\WW_0\otimes
\WW_{1,1}\otimes\cdots\otimes\WW_{1,p-1}\otimes\WW_{1,p+1}
\otimes\cdots\otimes\WW_{1,\ell_5}\mbox{ \ for \ }p\in\ol{1,\ell_5}.
\end{eqnarray}
In the following, an element $v\in\WW$ will always mean a homogeneous element with
$g(v)\in\Z_2$.
Observe from (\ref{Prod}) that
\begin{equation}
[u,vw]=[u,v]w+(-1)^{g(u)g(v)}v[u,w]\mbox{ \ for \ }u,v,w\in\WW.
\end{equation}
Thus for $a,b\in\WW_0,\,u,v\in\WW_1$ and $y,z\in\WW^\wedge_p$, we have
\begin{eqnarray}
\label{3.05-}
&\!\!\!\!\!\!\!\!&\!\!\!\!\!\!\!\!
[au,bv]=[a,b]uv+ba[u,v],\\
\label{3.05}
&\!\!\!\!\!\!\!\!&\!\!\!\!\!\!\!\![ys_p,zs_p]=[y\sptl_p,z\sptl_p]=0,
\\\label{3.05+}
&\!\!\!\!\!\!\!\!&\!\!\!\!\!\!\!\!
[ys_p\sptl_p,zs_p\sptl_p]=[ys_p\sptl_p,z]=[y,z]s_p\sptl_p,
\\
\label{3.06}
&\!\!\!\!\!\!\!\!&\!\!\!\!\!\!\!\!
[ys_p,z\sptl_p]=(-1)^{g(z)}\biggl([y,z]s_p\sptl_p+
(-1)^{g(y)g(z)}zy\biggr).
\end{eqnarray}

Now suppose  $\varphi$ is a $2$-cocycle satisfying
(\ref{2.6})--(\ref{2.8}).

{\bf Lemma 3.3.} {\it \
We have
\begin{eqnarray}
\label{3.07}
&\!\!\!\!\!\!\!\!\!\!&\!\!\!\!\!\!\!\!\!\!
\varphi(y,z)=\varphi(ys_p,z)=\varphi(ys_p,zs_p)=\varphi(y\sptl_p,z)
\nonumber\\
&\!\!\!\!\!\!\!\!\!\!&\ \ \ \
=
\varphi(y\sptl_p,z\sptl_p)=\varphi(ys_p\sptl_p,zs_p)=\varphi(ys_p\sptl_p,z\sptl_p)=0,
\\
\label{3.08}
&\!\!\!\!\!\!\!\!\!\!&\!\!\!\!\!\!\!\!\!\!
(-1)^{g(z)}\varphi(ys_p,z\sptl_p)=\varphi(y,zs_p\sptl_p)=\varphi(ys_p\sptl_p,zs_p\sptl_p),
\end{eqnarray}
for $p\in\ol{1,\ell_5}$ and $y,z\in\WW^\wedge_p$.
}

{\it Proof.} \
We have $\varphi(y,z)=\varphi([\sptl_p,s_py],z)=(-1)^{1+g(y)}\varphi(s_py,[\sptl_p,z])=0$.
The other equalities of (\ref{3.07}) follow from the fact
$$ 
0=\varphi(s_p\sptl_p,[ys_p^k\sptl_p^\mu,zs_p^{k'}\sptl_p^{\nu}])=
(k+k'-\mu-\nu)\varphi(ys_p^k\sptl_p^\mu,zs_p^{k'}\sptl_p^\nu)
\mbox{ \ (by (\ref{2.4}) and (\ref{2.8}))}
$$ 
for $k,k',\mu,\nu\in\{0,1\}$.
By (\ref{2.4}) and (\ref{2.7}), we have
$$
0=\varphi(\sptl_p,[ys_p,zs_p\sptl_p])=
(-1)^{g(y)}\biggl(\varphi(y,zs_p\sptl_p)-(-1)^{g(z)}\varphi(ys_p,z\sptl_p)\biggr),
$$
which gives the first equality of (\ref{3.08}).
Since $\WW^\wedge_p=[\WW^\wedge_p,
\WW^\wedge_p],$ by linearity, we can suppose $z=[z_1,z_2]$ for some $z_1,z_2\in\WW_p^\wedge$
without loss of generality.
Then by (\ref{3.05+}) and (\ref{2.4}),
$$
\begin{array}{ll}
\varphi(y,zs_p\sptl_p)\!\!\!\!&=
\varphi([y,z_1s_p\sptl_p],z_2s_p\sptl_p)
+(-1)^{g(y)g(z_1)}\varphi(z_1s_p\sptl_p,[y,z_2s_p\sptl_p])
\vs{4pt}\\&
=
\varphi([ys_p\sptl_p,z_1s_p\sptl_p],z_2s_p\sptl_p)
+(-1)^{g(y)g(z_1)}\varphi(z_1s_p\sptl_p,[ys_p\sptl_p,z_2s_p\sptl_p])
\vs{4pt}\\&
=\varphi(ys_p\sptl_p,[z_1s_p\sptl_p,,z_2s_p\sptl_p])
=\varphi(ys_p\sptl_p,zs_p\sptl_p),
\end{array}
$$
which gives the second equality of (\ref{3.08}).
\hfill$\Box$

Note that we have
\begin{equation}
\label{3.de}
\WW_1=[\WW_1,\WW_1]\oplus\F u_1,\mbox{ \ where \ }
u_1=s_1\cdots s_{\ell_5}\sptl_1\cdots\sptl_{\ell_5}\in\WW_1.
\end{equation}
We define a linear function $P:\WW_1\to\F$ by setting
\begin{equation}
\label{3.de1}
\mbox{
$P(u_1)=1$ and $P([\WW_1,\WW_1])=0$.
}
\end{equation}
Then by the second equation of (\ref{3.de1}), we have
\begin{equation}
\label{3.de1+}
P(uv)=(-1)^{g(u)g(v)}P(vu)\mbox{ \ \ for \ \ }u,v\in\WW_1.
\end{equation}
Define a bilinear function $\phi:\WW_0\times\WW_0\to\F$ by
setting
\begin{equation}
\label{3.de2}
\phi(a,b)=\varphi(au_1,bu_1)=\varphi(a,bu_1)\mbox{ \ \ for \ \ }a,b\in\WW_0.
\end{equation}
where the second equality follows from Lemma 3.3 and the fact that $u_1u_1=u_1$.
We have

{\bf Lemma 3.4.} {\it \ $\phi$ is a $2$-cocycle on $\WW_0$ satisfying $(\ref{2.6})$ and
$(\ref{2.7})$ $($with $\ell$ replaced by $\ell'_4$$)$. In particular $\phi=c\phi_0$ for some $c\in\F$ if $\ell'_4=\ell_3=1$, or $\phi=\sum_{\g\in\G}c_\g\phi_\g$ for some $c_\g\in\F$ if $\ell'_4=\ell_3=1$, or $\phi=0$ if $\ell'_2\ge1$ or $\ell'_4\ge2$.
}

{\it Proof.} \ The first statement can be verified directly, the second  follows from [S2] (we remark that although $\sum_{\g\in\G}c_\g\phi_\g$ may be an
infinite sum, it is {\it summable} in the sense that when it applies to any $(a,b)$
for $a,b\in\WW_0$, there are only finite many nonzero terms).
\hfill$\Box$

Our main result of this paper is the following. 

{\bf Theorem 3.5.} \ {\it
{\rm(1)} Suppose $\ell'_4=\ell_4=1$. Then $H^2(\WW,\F)= \F\ol\varphi_0$, where $\varphi_0$ is defined by
\begin{equation}
\label{Th3.4.1}
\varphi_0(au,bv)=\phi_0(a,b)P(uv)\mbox{ \ \ for \ \ }a,b\in\WW_0,\,u,v\in\WW_1,
\end{equation}
and
$\phi_0$ is defined by (\ref{Lem3.1}).

{\rm(2)} Suppose $ \ell'_4=\ell_3=1$. Then for any $ \g\in\G$, there corresponds a cohomology class $\ol\varphi_{\g}\in H^2(\WW,\F)$ defined by
\begin{equation}
\label{Th3.4.2}
\varphi_{\g}(au,bv)=\phi_\g(a,b)P(uv)\mbox{ \ \ for \ \ }
a,b\in\WW_0,\,u,v\in\WW_1,
\end{equation}
and $\phi_\g$ is defined by (\ref{Lem3.1.2}).
Furthermore, $H^2(\WW,\F)$ is a direct \vs{-5pt}product:
\begin{equation}\label{Theorom3.1.1}
H^2(\WW,\F)=\mbox{$\prod\limits_{\g\in\G}$}\F\ol\varphi_\g.\vs{-5pt}
\end{equation}

{\rm(3)} If $\ell'_2\ge1$ or $\ell'_4\ge2$, then $H^2 (\WW, \F) =0$.

Proof.} \
The result follows from Theorem 3.1 if $\ell_5=0$. Thus assume $\ell_5\ge1$.
First we verify that in case of $\ell'_4=\ell_4=1$,
a bilinear function $\varphi_0$ defined by (\ref{Th3.4.1})
is a nontrivial $2$-cocycle on $\WW$: The super-skew-symmetry follows from
the skew-symmetry of $\phi_0$ and (\ref{3.de1+}). By (\ref{3.05-}), we have
\begin{eqnarray}
\label{add-equ1}
\varphi_0(au,[bv,cw])&\!\!\!=\!\!\!&
\varphi_0(au,[b,c]vw+cb[v,w])
\nonumber\\&\!\!\!=\!\!\!&
\phi_0(a,[b,c])P(uvw)
+\phi_0(a,cb)P(u[v,w]),
\end{eqnarray}
for $a,b,c\in\WW_0,\,u,v,w\in\WW_1.$ Now the super-Jacobi identity follows from
(\ref{Rmk3.2}), (\ref{3.de1+}) and (\ref{add-equ1}) (together with its shifted version).
Clearly, $\varphi_0$ is nontrivial since $\phi_0$ is nontrivial.
Similarly, in case of $ \ell'_4=\ell_3=1$,
since (\ref{add-equ1}) still holds with $\varphi_0$ replaced by
$\varphi_\g$ for any $\g\in\G$,
one can prove that
$\varphi=\mbox{$\sum_{\g\in\G}$}
c_\g\varphi_\g$
is a nontrivial $2$-cocycle on $\WW$,
where $\varphi_\g$ is defined in
(\ref{Th3.4.2}), and $c_\g\in\F$ for $\g\in\G$ such that $c_\g\ne0$ for
at least one $\g\in\G$. In particular, the right-hand side of (\ref{Theorom3.1.1}) is a direct product.

Now suppose $\varphi$ is a $2$-cocycle satisfying (\ref{2.6})--(\ref{2.8}).
For any
\begin{equation}\label{u}
u=s_1^{j_1}\cdots s_{\ell_5}^{j_{\ell_5}}\sptl_1^{\mu_1}\cdots\sptl_{\ell_5}^{\mu_{\ell_5}}\in\WW_1,
\mbox{ \ \ where \ }j_p,\mu_p\in\{0,1\},
\end{equation}
we define its {\it support} to be
${\rm supp}(u)=\{p\!\in\ol{1,\ell_5}\,|\,(j_p,\mu_p)\!\ne\!0\}$.
We want to prove
\begin{eqnarray}
&&\varphi(a,bu)=0\mbox{ \ for \ }a,b\in\WW_0,\,u\in[\WW_1,\WW_1].
\label{End1}
\end{eqnarray}
Suppose $u\in[\WW_1,\WW_1]$ is as in (\ref{u}), then $u\ne u_1$ by (\ref{3.de}).
Thus there exists $p\in\ol{1,\ell_5}$ such that
$u=u'u''$ for some $u'\in\WW_p^\wedge$ and $u''=1,s_p$ or $\sptl_p$.
Now take $y=a,\,z=bu'\in\WW_p^\wedge$, by  (\ref{3.07}), we obtain (\ref{End1}).
Let $\phi$ be defined as in (\ref{3.de2}).
By using
(\ref{3.08}), (\ref{3.de1}), (\ref{3.de2}), (\ref{End1}) and induction on the support size $\#{\rm supp}(u)$,
one can similarly prove
\begin{eqnarray}
\label{End2}
&&\varphi(au,bv)=\varphi(a,buv)=\phi(a,b)P(uv)
\mbox{ \ for \ }a,b\in\WW_0,\,u,v\in\WW_1.
\end{eqnarray}
For example, if $\#{\rm supp}(u)=0$ then (\ref{End2}) follows from
(\ref{3.de1}), (\ref{3.de2}) and (\ref{End1}).

Now the theorem follows from (\ref{End2}) and Lemma 3.4.
\hfill$\Box$

\vskip9pt\ni{\bf References}\vskip5pt\small
\parindent=8ex\parskip=2pt\baselineskip=2pt
\re{BKV} B.~Bakalov, V.G.~Kac, A.A.~Voronov, ``Cohomology of
conformal algebras,'' {\it Comm.~Math. Phys.} {\bf200} (1999), 561--598.

\re{K1}  V.G.~Kac, {\it Infinite Dimensional Lie Algebras}, 3rd ed.,
Combridge Univ.~Press, 1990.
\re{K2} V.G.~Kac, ``Lie superalgebras,'' {\it Adv. Math.} {\bf26} (1977), 8-96.
\re{KP}  V.G.~Kac, D.H.~Peterson, ``Spin and pepresentation of
infinite dimensional Lie algebras and groups,'' {\it Proc.~Nat.~Acad.~Sci.
U.S.A.} {\bf78} (1981), 3308--3312.
\re{L}   W.~Li, ``2-Cocycles on the algebra of differential
operators,'' {\it J.~Algebra} {\bf122} (1989), 64--80.

\re{LW}   W.~Li, R.L.~Wilson, ``Central extensions of some Lie
algebras,'' {\it Proc.~Amer.~Math.~Soc.} {\bf126} (1998), 2569--2577.

\re{ScZ1} M.~Scheunert, R.B.~Zhang, ``Cohomology of Lie
  superalgebras and their generalizations,'' {\it J.~Math. Phys.} {\bf39}
  (1998), 5024--5061.

\re{ScZ2} M.~Scheunert, R.B.~Zhang, ``The second cohomology of $sl(m/1)$ with coefficients in its enveloping algebra is trivial,'' {\it Lett.~Math. Phys.} {\bf47}
  (1999), 33--48.

\re{S1} Y.~Su, ``2-Cocycles on the Lie algebras of all differential
operators of several indeterminates,'' {\it (Chinese) Northeastern Math.~J.} {\bf6} (1990), 365--368.

\re{S2} Y.~Su, ``2-cocycles on the Lie algebras of generalized
differential operators'', {\it Comm.~Algebra} {\bf30} (2002), 763--782.

\re{S3} Y.~Su, ``Classification of quasifinite modules over the Lie algebras
  of Weyl type,'' {\it Adv. Math.} {\bf174} (2003), 57--68.

\re{S4} Y.~Su, ``Low dimensional cohomology of general conformal algebras
$gc_N$,'' {\it J.~Math.~Phys.} {\bf45} (2004), 509--524.

\re{SXZ} Y.~Su, X.~Xu and H.~Zhang, ``Derivation-simple algebras
  and the structures of Lie algebras of Witt type,''
  {\it J.~Algebra} {\bf2000} (233), 642--662.

\re{SZ1} Y.~Su, K.~Zhao, ``Simple Lie algebras of Weyl type,'' {\it Science in China
A} {\bf 44} (2001), 419--426.

\re{SZ2} Y.~Su, K.~Zhao, ``Isomorphism classes and automorphism
  groups of algebras of Weyl type,'' {\it Science in China A} {\bf 45} (2002),
  953--963.

\re{SZ3}  Y.~Su, K.~Zhao, ``Second cohomology group of generalized Witt
type Lie algebras and certain reperesentations,'' {\it Comm.~Alegrba} {\bf30} (2002),  3285--3309.

\re{SZZ} Y.~Su, K.~Zhao and L.~Zhu,
    ``Simple Lie color algebras of Weyl types,'' {\it Israel J. Math.}
    in press.
\re{X1} X.~Xu, ``New generalized simple Lie algebras of Cartan
    type over a field with characteristic 0,''
   {\it J.~Algebra} {\bf 224} (2000), 23--58.
\re{X2} X.~Xu, ``Generalizations of Block algebras,''
   {\it Manuscripta Math.} {\bf100} (1999), 489--518.
\re{X3} X.~Xu,  ``Simple conformal superalgebras of finite
  growth,'' {\it Algebra Colloquium} {\bf7} (2000), 205--240.
\re{X4} X.~Xu, ``Equivalence of conformal superalgebras to
   Hamiltonian superoperators,'' {\it Algebra Colloquium} {\bf8} (2001), 63--92.
\re{X5} X.~Xu, ``Quadratic Conformal Superalgebras,'' {\it
    J.~Algebra} {\bf231} (2000), 1--38.
\re{Z} K.~Zhao,
    ``Simple algebras of Weyl type II,"
   {\it Proc. Amer. Math. Soc.} {\bf130} (2002), 1323--1332.
\end{document}